\documentclass[11pt]{article}
\usepackage{amsthm}
\usepackage{fullpage}
\usepackage{graphicx}
\usepackage{amsmath}
\usepackage{pstricks}
\newpsobject{showgrid}{psgrid}{subgriddiv=1,griddots=10,gridlabels=6pt}
\newcommand{\Chi}{\mbox{\large $\mathbf{\chi}$}}

\newcommand {\abs} [1] {\left| #1 \right|}

\newcommand {\lip}{\mathop{\mathrm{Lip}}}

\newtheorem {thm} {Theorem}[section]

\newtheorem {lem} [thm] {Lemma}

\numberwithin{equation}{section}

\title{Stable Manifolds and Homoclinic Points Near Resonances in the Restricted Three-Body Problem}
\author{D. Viswanath \thanks{Department of Mathematics, University of Michigan, 
530 Church Street, Ann Arbor, MI 48109, U.S.A. This work was supported by the NSF
grant DMS-0407110
and by a research fellowship from the Sloan Foundation.}}

\begin{document}
\maketitle
\begin{abstract}
The restricted three-body problem describes the motion of a massless particle under the influence of two primaries of masses $1-\mu$ and $\mu$ that circle each other with period
equal to $2\pi$. For small $\mu$, a resonant periodic motion of the massless particle in 
the rotating frame can be described by relatively prime integers $p$ and $q$, if its 
period around the heavier primary is approximately $2\pi p/q$, and by its approximate eccentricity $e$. We give a method for the formal development of the stable and unstable
manifolds associated with these resonant motions. We prove the validity of this formal 
development and the existence of homoclinic points in the resonant region.

In the study of the Kirkwood gaps in the asteroid belt, the separatrices of the averaged
equations of the restricted three-body problem are commonly used to derive analytical
approximations to the boundaries of the resonances. We use the unaveraged equations to find values of asteroid eccentricity below which these approximations
will not hold for the Kirkwood gaps with $q/p$ equal to $2/1$, $7/3$, $5/2$, $3/1$,
and $4/1$. 

Another application is to the existence of asymmetric librations in the exterior resonances.
We give values of asteroid eccentricity below which asymmetric librations will not exist
for the $1/7$, $1/6$, $1/5$, $1/4$, $1/3$, and $1/2$ resonances for any $\mu$ however small.
But if the eccentricity exceeds these thresholds, asymmetric librations will exist for 
$\mu$ small enough in the unaveraged restricted three-body problem.

\end{abstract}
\section{Introduction}

The restricted three-body problem describes the motion of a massless particle under the
influence of two primaries of masses $1-\mu$ and $\mu$. The Hamiltonian
\begin{equation}
H = \frac{1}{2}(p_x^2+p_y^2)+yp_x - xp_y - \frac{1-\mu}{(x^2+y^2)^{1/2}}
-\mu\Biggl(\frac{1}{((x-1)^2+y^2)^{1/2}}-x\Biggr)
\label{eqn-1-1}
\end{equation}
gives the equations of motion of the massless particle.
In \eqref{eqn-1-1}, it is assumed that the primary of mass $1-\mu$ is located at
$(0,0)$, that the primary of mass $\mu$ is located at $(1,0)$, and that
the frame of reference rotates with the second primary 
in the anticlockwise sense with period $2\pi$.

In terms of the heliocentric Delaunay variables $L,l,G,g$, the Hamiltonian becomes
\begin{equation}
H = -\frac{(1-\mu)^2}{2L^2}-G-\mu\Omega'(L,G,l,g)
= -\frac{1}{2L^2}-G-\mu\Omega(L,G,l,g) -\frac{\mu^2}{2L^2},
\label{eqn-1-2}
\end{equation}
where
\begin{align}
\Omega' &= \frac{1}{(1+r^2-2r\cos\theta)^{1/2}}-r\cos\theta\nonumber\\
\Omega &= \Omega'-1/L^2.
\label{eqn-1-3}
\end{align}
The variables $r,\theta$ can be obtained in terms of $L,l,G,g$ using the equations
\begin{align}
e &= (1-G^2/L^2)^{1/2}\nonumber\\
a &= L^2/(1-\mu)\nonumber\\
l &= E-e\sin E\nonumber\\
\cos \nu &= (\cos E-e)/(1-e\cos E)\nonumber\\
\sin \nu &= (1-e^2)^{1/2}\sin E/(1-e\cos E)\nonumber\\
\theta &= g + \nu \nonumber\\
r &= a(1-e\cos E).
\label{eqn-1-4}
\end{align}
The Cartesian coordinates $x,y$ in the rotating frame used in \eqref{eqn-1-1} are given by
$x=r\cos\theta$ and $y=r\sin\theta$. When $\mu=0$, the orbit of the massless particle in the inertial frame is an ellipse with eccentricity $e$ and semimajor axis $a$;
$l$, $\nu$, and $E$ are the mean, true, and eccentric anomalies; $g$ is the argument of the perihelion in the rotating frame and $G$ denotes angular momentum. The variables $r$ and $\theta$ are the polar coordinates of the massless particle in the rotating frame.

We use $L,l,G,g$ to investigate motion near resonances as the Hamilton's equations
of \eqref{eqn-1-2} take an especially simple form for $\mu$ small. If $L=(p/q)^{1/3}$,
$G=(p/q)^{1/3}(1-e^2)^{1/2}$, $0<e<1$,
and $\mu=0$, the massless particle moves on an ellipse with period $2\pi p/q$ and eccentricity $e$ in the inertial frame. It is assumed that $p$ and $q$ are relatively prime positive integers. The motion is periodic in the rotating frame as well with period equal to $2\pi p$. If $l$ and $g$ are multiples of $\pi$ and the unperturbed orbit does not collide with the second primary, these periodic motions in the rotating frame perturb to periodic motions for $\mu > 0$ and $\mu$ small \cite{Barrar1}. Given
$p,q,e$ it might appear that there are four possibilities as $l$ and $g$ can be either
$0$ or $\pi$, but in fact only two of these are distinct. These are the two  $q/p$
resonant periodic motions of the restricted three-body problem that correspond to 
eccentricity $e$.

\begin{figure}
\begin{center}
\psset{xunit=1.3cm,yunit=0.45cm}
\begin{pspicture}(-1,-1)(7,10)
\psline[linewidth=0.5mm](0,0)(6,0)
\psline[linewidth=0.5mm](0,0)(0,10)
\psline[linewidth=0.5mm](6,0)(6,10)
\rput(3,-0.5){\large $l$}
\rput(-0.3,6){\large $L$}
\rput(0,-0.5){$0$}
\rput(6,-0.5){$2\pi$}
\psline[linestyle=dotted](0,3)(6,3)
\psline[linestyle=dotted](0,5)(6,5)
\psline[linestyle=dotted](0,7)(6,7)
\psline[linestyle=dotted](0,9)(6,9)
\rput(-0.5,4){$\bigl(\frac{2}{7}\bigr)^{1/3}$}
\rput(-0.5,8){$\Bigl(\frac{1}{3}\Bigr)^{1/3}$}
\qdisk(3,4){0.03}
\pscircle(3,4){0.15}
\psarc[fillcolor=black](0,4){0.03}{-90}{90}
\psarc(0,4){0.15}{-90}{90}
\psarcn[fillcolor=black](6,4){0.03}{-90}{90}
\psarcn(6,4){0.15}{-90}{90}
\qline(1.35,4.15)(1.65,3.85)
\qline(1.35,3.85)(1.65,4.15)
\qline(4.35,4.15)(4.65,3.85)
\qline(4.35,3.85)(4.65,4.15)
\pscurve(4.5,4)(5.25,4.5)(6,4.9)
\pscurve[linestyle=dashed](0,4.9)(0.75,4.5)(1.5,4)
\qdisk(3,8){0.03}
\pscircle(3,8){0.15}
\qline(0,8)(0.15,8.15)
\qline(0,8)(0.15,7.85)
\qline(6,8)(5.85,8.15)
\qline(6,8)(5.85,7.85)
\pscurve(0,8)(1.5,8.7)(3,8.9)
\pscurve[linestyle=dashed](6,8)(4.5,8.7)(3,8.9)
\end{pspicture}
\psset{xunit=1cm,yunit=1cm}
\end{center}
\label{fig-1}
\caption[xyz]{The $l$-$L$ plane shown above is the Poincar\'e section $g=0$.
Periodic points as well as stable and unstable manifolds are shown for
$q/p = 7/2$ and $q/p=3/1$.}
\end{figure}
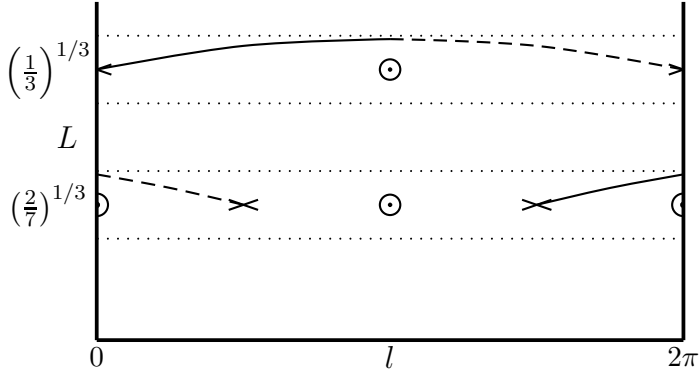

The Hamiltonian $H$ is conserved by the flow, and for a given $H$ and $\mu$ small it is possible to solve \eqref{eqn-1-2} for $G$ using the implicit function theorem. 
Thus the Poincar\'e section $g=0$ can be identified with the $l$-$L$ plane. Each resonant periodic motion corresponding to $p,q,e$ intersects this Poincar\'e section exactly $p$ times. Typically, if one of the resonant periodic motions is of elliptic type the other is of hyperbolic type \cite{Viswanath1}. In Figure \ref{fig-1}, we have shown the periodic points on the $L$-$l$ section  for $q/p =3/1$ and $q/p=7/2$ with elliptic points marked as circles and hyperbolic points as crosses. The perturbing term $\Omega$ and
the $O(\mu^2)$ term in \eqref{eqn-1-2} are unchanged by the transformation $L\leftarrow L$, $l\leftarrow -l$, $G\leftarrow G$, $g\leftarrow -g$, which is symplectic with multiplier $-1$. This discrete symmetry of the Hamiltonian \eqref{eqn-1-2} has the following implication for the return map to the $l$-$L$ section: if $(l_0,L_0)$ maps to $(l_1, L_1)$ then the return map sends $(-l_1,L_1)$ to $(-l_0, L_0)$.

We obtain a scaled version of the return map near $q/p$ resonance in Sections 2 and 3. 
The stable and unstable manifolds of the hyperbolic points of such a return map nearly
coincide and the angle of transversality can be upper bounded by a quantity that is
exponentially small in the small parameter $\mu$. A discussion of exponential splitting
of separatrices can be found in the work of Gelfreich and Lazutkin \cite{GL1},
Fontich and Sim\'o \cite{FS1}, and Holmes, Marsden and Scheurle \cite{HMS1}. Instability is
often associated with resonance, and it is therefore natural to look for transverse 
homoclinic points near resonances. However, it appears that there has been no construction
of transverse homoclinic points near resonances in either the restricted three-body problem or in some other
version of the planetary problem.
The rescaling of the return map given in Section 3 brings this 
problem into sharper focus and it is possible to make an analogy to the discussion of the
standard map in \cite{GL1}.

In Section 4, we describe a procedure for the formal development of the stable manifold
of a resonant periodic solution in powers of $\mu^{1/2}$. This procedure is specially 
adapted to the restricted three-body problem. In Section 5, we have included a proof of
validity of this formal expansion. We also prove the existence of homoclinic points.
A part of the verification essential for the proof is carried out in Section 6.

Although  investigations of the Kirkwood gaps in the asteroid belt have used physical
models that include the secular variation of Jupiter's elements and the effect of Saturn,
the averaged circular restricted three-body problem is still used to approximate the
boundaries of the resonant regions. The more complicated models are essential to explain the
dynamics within the resonant regions. An account of these models and their use can be
found in the monograph by Morbidelli \cite{Morbidelli1}. 
The use of the averaged equations of the circular restricted problem to sketch the boundaries
of resonance can be found in Dermott and Murray \cite{DM1}, Henrard and Lema\^itre \cite{HL1},
and Lema\^itre \cite{Lemaitre1}.
It appears to be known that the boundaries
obtained from the averaged circular restricted problem do not work well at low eccentricities
\cite{Morbidelli1} \cite{Yoshikawa1}. In fact for the $q/p=2/1$ case, the left boundary cannot
even be computed near $e=0$ \cite{Morbidelli1}. In Section 7, we  interpret the calculation of the boundaries
in terms of the unaveraged circular restricted problem. For the commonly studied Kirkwood 
gaps, we give values of the eccentricity below which the approximation of the boundary will not be valid.

Asymmetric librations near exterior resonances with $p > q$
in the averaged circular restricted problem have been investigated
by Beaug\'e \cite{Beauge1} and more recently by Voyatis et al. \cite{VKH1}. In Section 8,
we show the existence of these librations in the unaveraged circular restricted problem
for $q/p$ equal to $1/2$, $1/3$, $1/4$, $1/5$, $1/6$, and $1/7$. We also give minimum 
values of eccentricity for each of these exterior resonances which must be exceeded for
the asymmetric librations to exist. 

\section{Perturbative form of the return map}
 
Let $L$ range over the interval $[(p/q)^{1/3}-\delta_L,(p/q)^{1/3}+\delta_L]$ for some $\delta_L > 0$. We assume $p/q\neq 1/1$.
Let $e$ range over the interval $[e_{min}, e_{max}]$ such that $L^2 (1-e) > 1+\delta$ if $p>q$ and such that $L^2(1+e) < 1-\delta$ if $p<q$, for all allowed values of $L$ and $e$ and some $\delta > 0$. 
In addition, assume $0<e_{min}$ and $e_{max}<1$. 
We take the range of the Hamiltonian $H$ to be the set of values of $-(p/q)^{-2/3}/2-(p/q)^{1/3}(1-e^2)$ for $e_{min}\leq e\leq e_{max}$.
Then for any allowed value of $H$ and $L$, any real values of the angles $l$ and $g$ measured modulo $2\pi$, and $\mu$ sufficiently small, the implicit function theorem enables us to solve \eqref{eqn-1-2} for $G$ uniquely. In fact, $G$ will be an analytic function of $l,L,g,H,\mu$.

We use either $g=0$ or $g=\pi$ to define the Poincar\'e section. Since Hamilton's equations formed using the Hamiltonian \eqref{eqn-1-2} imply that $dg/dt = -1 + O(\mu)$, the return map is well defined for $L$ and $H$ in the intervals specified by the previous paragraph, any value of $l$, and $\mu$ small. Since $H$ is conserved by the flow, we may identify the Poincar\'e section for fixed $H$ with the $l$-$L$ plane as in Figure \ref{fig-1}. The return map preserves the area element $dl\,dL$.  The return map will be denoted by $T_1$ and $T_1^p$ will be denoted by $T_{p}$. In this section, we will obtain the perturbative form of $T_{p}$.

The Hamilton's equations of \eqref{eqn-1-2} imply that $\dot{l} = 1/L^3-\mu \Omega_L + O(\mu^2)$,
$\dot{g} = -1 - \mu \Omega_G$, and $\dot{L} = \mu \Omega_l$. We seek a solution of these equations with the initial conditions $l(0)=l_0$, $L(0) = L_0$, $g(0)=g_0$, where
$g_0$ is either $0$ or $\pi$, and $G(0) = G_0$. It is understood that $G$ must be obtained by solving \eqref{eqn-1-2}. If the solution is represented as $l(t) = l_a(t)
+l_b(t) \mu + O(\mu^2)$, $g(t) = g_a(t) +  g_b(t)\mu + O(\mu^2)$, and $L(t) = L_a(t) + L_b(t)\mu + O(\mu^2)$, then $l_a(t) = l_0+t/L_0^3$, $g_a(t) = g_0 - t$, and
$L_a(t) = L_0$. In addition, $l_b$, $g_b$, and $L_b$ must satisfy
\begin{equation*}
\dot{l}_b = (-3/L_a^4)L_b - \Omega_L,\quad \dot{g}_b = -\Omega_G,
\quad \dot{L}_b = \Omega_l,
\end{equation*}
where the partial derivatives of  $\Omega$ must be evaluated at $l=l_0+t/L_0^3$, $g=g_0-t$, $L = L_0$, and $G=G_0$. By solving the equations above, we get
\begin{align}
l(t) &= l_0 + t/L_0^3 + \mu\Biggl(-\int_0^t\Omega_L\, dt - 
\frac{3}{L_0^4} \int_0^t\int_0^\tau \Omega_l\, dt\, d\tau\Biggr) + O(\mu^2)\nonumber\\
g(t) &= g_0 - t -\mu \int_0^t \Omega_G\, dt + O(\mu^2)\nonumber\\
L(t) &= L_0 + \mu \int_0^t \Omega_l\, dt + O(\mu^2),
\label{eqn-2-1}     
\end{align}
where the partial derivatives of  $\Omega$ must be evaluated at $l=l_0+t/L_0^3$, $g=g_0-t$, $L = L_0$, and $G=G_0$. The solution given by \eqref{eqn-2-1} is valid over any finite interval of time for initial conditions in the domain already indicated and for $\mu$ sufficiently small.

To approximate $T_{p}$, it is necessary to find the time $t_r$ at which $g(t) = -2\pi p + g_0$. Using the equation for $g(t)$ in \eqref{eqn-2-1}, we get
\(t_r = 2\pi p - \mu \int_0^{2\pi p}\Omega_G\, dt + O(\mu^2).\)
Using \eqref{eqn-2-1}, we may deduce that $T_p$ is given by
\begin{align}
l_1 &= l_0 + 2\pi p/L_0^3 + \mu \Biggl(-\frac{1}{L_0^3}\int_0^{2\pi p} \Omega_G\, dt
- \int_0^{2\pi p} \Omega_L\, dt 
- \frac{3}{L_0^4} \int_0^{2\pi p} \int_0^\tau \Omega_l\,dt\,d\tau\Biggr) + O(\mu^2)
\nonumber\\
L_1 &= L_0 + \mu \int_0^{2\pi p} \Omega_l\, dt + O(\mu^2),
\label{eqn-2-2}
\end{align}
where the partial derivatives of  $\Omega$ must be evaluated at $l=l_0+t/L_0^3$, $g=g_0-t$, $L = L_0$, and $G=G_0$. The expression for $T_p$ given by \eqref{eqn-2-1} is valid for
$L \in [(p/q)^{1/3}-\delta_L, (p/q)^{1/3}+\delta_L]$, for any real value of the angle $l$, for $H$ within a range that ensures avoidance of collision with the second primary as specified earlier, and for $\mu$ sufficiently small. Since this domain of validity is compact, the $O(\mu^2)$ terms in \eqref{eqn-2-2} hold uniformly over the domain. 

The lemmas below are related to the first return map $T_1$ and its $p$th iterate $T_p$.

\begin{lem}
Assume that $(l_0,L_0)$ maps to $(l_1, L_1)$ under $T_1$ (or $T_p$). Then $(-l_1,L_1)$ maps to $(-l_0,L_0)$ under $T_1$ (or $T_p$).
\label{lem-2-1}
\end{lem} 
\begin{proof}
The Hamiltonian \eqref{eqn-1-2} is unchanged by the transformation $l\leftarrow -l$,
$L \leftarrow L$, $g\leftarrow -g$, $G\leftarrow G$. Therefore, if $l(t),\, g(t),\,L(t),\, G(t)$ is a solution of the Hamilton's equations of \eqref{eqn-1-2} for $0\leq t \leq t^\ast$, then $-l(-t),\, -g(-t),\, L(-t),\, G(-t)$, where 
$-t^\ast\leq t \leq 0$, is also a solution. The lemma follows if it is noted that the Poincar\'e section is defined using either $g=0$ or $g=\pi$.
\end{proof}

The lemma below is useful for finding fixed points of $T_p$.
\begin{lem} Let
\[ \phi_{p}(l_0, L_0) = \int_0^{2\pi p} \Omega_l\, dt,\]
where the arguments of $\Omega_l$ is evaluated with $l=l_0+t/L_0^3$, $g=g_0-t$, $L = L_0$, and $G=G_0$. Then $\phi_p(l_0, (p/q)^{1/3}) = \phi_p(l_0 + 2\pi/p, (p/q)^{1/3})$ and
$\phi_p(-l_0, (p/q)^{1/3}) = -\phi_p(l_0, (p/q)^{1/3})$.
\label{lem-2-2}
\end{lem}
\begin{proof}
If we define $\phi_1(l_0, L_0) = \int_0^{2\pi} \Omega_l\, dt$, then $L_1 = L_0 + 
\phi(l_0, L_0) \mu + O(\mu^2)$ under the first return map $T_1$. Using Lemma \ref{lem-2-1}, we may conclude that $\phi_1(l_0, L_0)+\phi(-l_0-2\pi/L_0^3,L_0) = 0$.
If $L_0 = (p/q)^{1/3}$, we have 
\begin{equation}
\phi_1(l_0, (p/q)^{1/3})+\phi_1(-l_0-2\pi q/p, (p/q)^{1/3}) = 0.
\label{eqn-2-3}
\end{equation}  
From $T_p = T_1^p$, we get
\begin{equation}
\phi_p(l_0, L_0) = \sum_{j=0}^{p-1}\phi_1(l_0 + 2j\pi/L_0^3).
\label{eqn-2-4}
\end{equation}
That $\phi_p(l_0, (p/q)^{1/3})$ has period equal to $2\pi/p$ in $l_0$ follows from
\eqref{eqn-2-4} and elementary number theory if it is noted that $p$ and $q$ are relatively prime. That $\phi_p(l_0, (p/q)^{1/3})$ is an odd function of $l_0$ follows from \eqref{eqn-2-3}. 
\end{proof}

\section{Scaling and periodic points of the return map}
The variable $\lambda$ defined by $L = (p/q)^{1/3} + \lambda \sqrt{\mu}$ can be used instead of $L$ to blow up the region of the Poincar\'e section near $p/q$ resonance. 
The image of the point $(l_0,\lambda_0)$ in  the $l$-$\lambda$ plane under the $p$th return map $T_p$ can be calculated using \eqref{eqn-2-2}, and it is given by
\begin{align}
l_1 &= l_0 - c_1 \lambda_0 \mu^{1/2} + (c_2 \lambda_0^2 + \chi(l_0))\mu + O(\mu^{3/2})
\nonumber\\
\lambda_1 &= \lambda_0 + \phi(l_0)\mu^{1/2} + \lambda_0\psi(l_0)\mu + O(\mu^{3/2}),
\label{eqn-3-1}
\end{align}
where
\begin{align}
c_1 &= 6\pi q^{4/3}p^{-1/3}\quad\text{and}\quad 
c_2 = 12\pi q^{5/3}p^{-2/3} \nonumber\\
\phi(l_0) &= \int_0^{2\pi p} \Omega_l\, dt\nonumber\\
\psi(l_0) &= \int_0^{2\pi p} \Omega_{lL}\, dt + 
\frac{q}{p} \int_0^{2\pi p}\Omega_{lG}\, dt 
-\frac{3q^{4/3}}{p^{4/3}}\int_0^{2\pi p} t \Omega_{ll}\, dt\nonumber\\
\chi(l_0) &= -\int_0^{2\pi p} \Omega_L\, dt
-\frac{q}{p}\int_0^{2\pi p}\Omega_G\, dt 
-\frac{3q^{4/3}}{p^{4/3}} \int_0^{2\pi p}\int_0^{\tau}\Omega_l\,dt\,d\tau.
\label{eqn-3-2}
\end{align}
In \eqref{eqn-3-2}, the partial derivatives of $\Omega$ must be evaluated at
$l=l_0 + qt/p$, $g = g_0-t$ (where $g_0$ is either $0$ or $\pi$ depending upon the choice of the Poincar\'e section), $L=(p/q)^{1/3}$, and $G = (p/q)^{1/3}(1-e^2)^{1/2}$.
The expression for $T_p$ given by \eqref{eqn-3-1} and \eqref{eqn-3-2} is valid for any real value of $l_0$ and $\abs{\lambda_0} \leq \delta_L/\mu^{1/2}$. The domain of definition can therefore be taken as $\abs{\lambda_0} \leq C_\lambda$ with any positive constant $C_\lambda$ for sufficiently small $\mu$. 
To derive the expression for $\psi(l_0)$ given in \eqref{eqn-3-2}, we must use the second line of \eqref{eqn-2-2} and notice that
the equation $L_0 = (p/q)^{1/3}+\lambda_0\mu^{1/2}$ and \eqref{eqn-1-2} imply $G_0 = (p/q)^{1/3}(1-e^2)^{1/2}+(q/p)\lambda_0\mu^{1/2} + O(\mu)$. A term equal to $2\pi q$ has been dropped from the first line of \eqref{eqn-3-1} as $l_0$ and $l_1$ are angles measured modulo $2\pi$.

The function $\phi(l_0)$ equals $\phi_p(l_0, L_0)$ defined by Lemma \ref{lem-2-1} when $L_0 = (p/q)^{1/3}$. By Lemma \ref{lem-2-1}, $\phi(l_0)$ is an odd function with period equal to $2\pi/p$. Therefore, $\phi(0)=0$ and $\phi(\pi/p) = 0$. We make the following 
assumption about $\phi(l_0)$:
\begin{quote}
Assumption A: For $l_0\in [0,2\pi/p)$, the only points where $\phi(l_0)=0$ are $l_0 = 0$ and $l_0=\pi/p$. At these points, the derivative $\phi'(l_0)$ is nonzero.	
\end{quote}

We turn to the verification of this assumption in Section 6. An example of a function which is odd with period $2\pi/p$ and which satisfies the assumption above is $\sin(p l_0)$. In fact, it will be shown later that $\phi(l_0)$ is proportional to $\sin(p l_0) e^{\abs{p-q}}$ for small $e$.

The assumption about $\phi(l_0)$ can be put to use to find fixed points of $T_p$ in the $l$-$\lambda$ plane. We can use \eqref{eqn-3-1} and write
\begin{align}
(l_1-l_0)/\mu^{1/2} &= -c_1\lambda_0 + (c_2\lambda_0^2 + \chi(l_0))\mu^{1/2} + O(\mu)
\nonumber\\
(\lambda_1-\lambda_0)/\mu^{1/2} &= \phi(l_0) + \lambda_0\psi(l_0)\mu^{1/2} + O(\mu).
\label{eqn-3-3}
\end{align}
When $\mu=0$, the right hand sides of the two equations in \eqref{eqn-3-3} are both zero 
if $\lambda_0 = 0$ and if $l_0$ is an integral multiple of $\pi/p$. The implicit function theorem, along with the assumption about $\phi(l_0)$ stated above, allows us to infer that
the right hand sides in \eqref{eqn-3-3} are $0$ for $\mu$ sufficiently small, 
\begin{equation}
l_0 = j\pi /p + O(\mu),\quad \text{and}\quad 
\lambda_0 = \bigl(\chi(j\pi /p)/c_1\bigr)\mu^{1/2} + O(\mu),
\label{eqn-3-4} 
\end{equation}
where $j$ is an integer. Thus the points given by \eqref{eqn-3-4} are fixed points of $T_p$ for $\mu$ sufficiently small. Since each application of the first return map $T_2$ increments $l_0$ by $2\pi q/p + O(\mu)$, where $q$ and $p$ are relatively prime, we may group the fixed points given by \eqref{eqn-3-4} into two sets, the first with 
$j=0,2,\ldots, 2(p-1)$ and the second with $j=1,3,\ldots,2(p-1)+1$. Then any fixed point in the first set moves to all other points in that set upon successive applications of $T_2$ and returns to itself after the $p$th application; and likewise with the second set. The assumption about $\phi(l_0)$ implies that if one set of periodic points is elliptic then the other set is hyperbolic as will become clear shortly.

If Assumption A about $\phi(l)$ fails to hold, points $l_0$ with 
$\phi(l_0) = 0$ and $\phi'(l_0) \neq 0$ will still correspond to
periodic points of \eqref{eqn-3-1} for $\mu$ sufficiently small. 
If $l_0\neq 0$ and $l_0 \neq \pi$, the corresponding periodic points can
be the centers of asymmetric librations as described in the last section.

\section{Formal expansion of the stable manifold}
The expression for the map $T_p$ given by \eqref{eqn-3-1} and \eqref{eqn-3-2} can be rewritten by shifting the center of the $l$-$\lambda$ plane to $(l^\ast, \lambda^\ast)$,
where $(l^\ast, \lambda^\ast)$ is the fixed point of $T_p$ given by \eqref{eqn-3-4}, with $j$ being some integer. The map $T_p$ applied to the $l$-$\lambda$ plane centered at such a fixed point takes the form
\begin{align}
l_1 &= l_0 - c_1 \lambda_0 \mu^{1/2} + (c_2 \lambda_0^2 + \Chi(l_0))\mu 
+ r(l_0, \lambda_0, \sqrt{\mu})\mu^{3/2} \nonumber\\
\lambda_1 &= \lambda_0 + \Phi(l_0)\mu^{1/2} + \lambda_0 \Psi(l_0) \mu 
+ s(l_0, \lambda_0, \sqrt{\mu})\mu^{3/2}, 
\label{eqn-4-1}
\end{align}
where
\begin{equation}
\Chi(l_0) = \chi(l_0+j\pi/p)-\chi(j\pi /p),\quad
\Phi(l_0) = \phi(l_0+j\pi/p),\quad
\Psi(l_0) = \psi(l_0+j\pi/p).
\label{eqn-4-2}
\end{equation}
The derivations of \eqref{eqn-2-2} and \eqref{eqn-3-1} imply that the remainder terms $r$ and $s$ in \eqref{eqn-4-1} are analytic in $l_0$, $\lambda_0$, and $\sqrt{\mu}$ for
any real $l_0$, $\lambda_0$  with $\abs{\lambda_0} \leq C_\lambda$, and $\sqrt{\mu}$ sufficiently small in magnitude. Besides, $r(0,0,\sqrt{\mu}) = s(0,0,\sqrt{\mu})=0$ because $(0,0)$ is a fixed point of \eqref{eqn-4-1}.

The Jacobian $dT_p$ of \eqref{eqn-4-1} at the origin is given by
\begin{equation}
dT_p = \begin{pmatrix} 1 & -c_1\mu^{1/2}\\\Phi'(0)\mu^{1/2} & 1 \end{pmatrix}
+ \begin{pmatrix} \Chi'(0) & \\ & \Psi(0)\end{pmatrix} \mu + O(\mu^{3/2}).
\label{eqn-4-3}
\end{equation}
For $\mu>0$ and $\mu$ small, the fixed point is hyperbolic if $\Phi'(0) < 0$ and elliptic if $\Phi'(0) > 0$. From \eqref{eqn-4-2}, it follows that $\Phi'(0) = \phi'(j\pi /p)$.
The assumption of Section 3 implies that $\phi'(0)$ and $\phi'(\pi/p)$ are of opposite signs. Thus if the set of fixed points given by \eqref{eqn-3-4} is of hyperbolic or elliptic type for even $j$, the set of fixed points given by odd $j$ must be of the opposite type. We shall assume that the fixed point used to shift the coordinate system and obtain \eqref{eqn-4-1} to be of hyperbolic type, which means $\Phi'(0) < 0$. Let $\alpha = \sqrt[4]{-\Phi'(0)/c_1}$. Then a calculation using \eqref{eqn-4-3} shows that the eigenvalues $1-\alpha^2 c_1\mu^{1/2} + \bigl((\Chi'(0)+\Psi(0))/2\bigr)\mu
+ O(\mu^{3/2})$ and $1+\alpha^2 c_1\mu^{1/2} + \bigl((\Chi'(0)+\Psi(0))/2\bigr)\mu
+ O(\mu^{3/2})$ of $dT_p$ correspond to eigenvectors of slopes
\begin{equation}
\alpha^2+\frac{\Chi'(0)-\Psi(0)}{2c_1} \mu^{1/2} + O(\mu)\quad \text{and}\quad
-\alpha^2+\frac{\Chi'(0)-\Psi(0)}{2c_1} \mu^{1/2} + O(\mu),
\label{eqn-4-4}
\end{equation}
respectively. The slope of the stable manifold of the fixed point of the map \eqref{eqn-4-1} at the origin must be given by the first of the two expressions in \eqref{eqn-4-4}. When we derive an approximation to that stable manifold, \eqref{eqn-4-4} will serve to check  the correctness of that approximation.

\begin{figure}
\begin{center}
\includegraphics[scale=0.7]{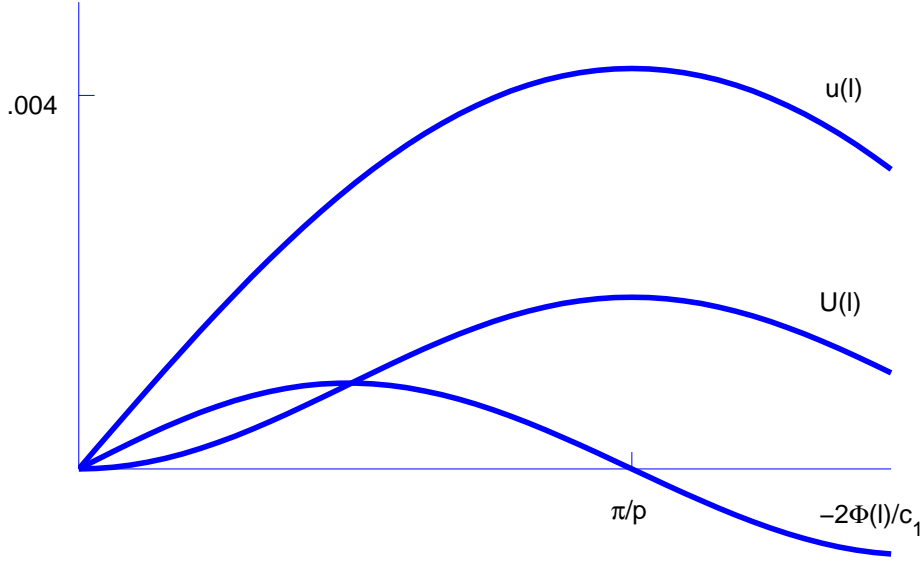}
\end{center}
\caption[xyz]{In the figure above, $u(l)$ has been scaled down by $0.1$ to make it fit. The plots correspond to the case $q/p=3/1$, $e=0.1$.}
\label{fig-2}
\end{figure}

To find the stable manifold of this fixed point, we rewrite
\eqref{eqn-4-1} in the following form:
\begin{align}
\lambda_0 \mu^{1/2} &=
-\frac{l_1-l_0}{c_1} + \frac{c_2\lambda_0^2+\Chi(l_0)}{c_1}\mu
+\cdots \nonumber\\
\lambda_1^2 &= \lambda_0^2 + 2\lambda_0\Phi(l_0)\mu^{1/2}
+ (\Phi(l_0)^2+2\lambda_0^2\Psi(l_0))\mu +\cdots
\label{eqn-4-5}
\end{align}
Let $(l_0,\lambda_0)$ be a point on the stable manifold. Its
iterates $(l_1, \lambda_1)$, $(l_2,\lambda_2)$, and so on are also
on the stable manifold. In addition, $l_n \rightarrow 0$
and $\lambda_n\rightarrow 0$ as $n\rightarrow\infty$. Using
\eqref{eqn-4-5}, we may write
\begin{align*}
\lambda_n^2 &= \lambda_0^2+ 2
\bigl( \lambda_0\Phi(l_0) + \lambda_1\Phi(l_1)+\cdots
+\lambda_{n-1}\Phi(l_{n-1})\bigr)\mu^{1/2} + O(\mu)\\
	&= \lambda_0^2-\frac{2}{c_1}
\bigl( (l_1-l_0)\Phi(l_0) + (l_2-l_1)\Phi(l_1)+\cdots
+(l_n-l_{n-1})\Phi(l_{n-1})\bigr) + O(\mu).
\end{align*}
Note that 
\begin{equation}
(l_{j+1}-l_j)\Phi(l_j) =
\int_{l_j}^{l_{j+1}} \Phi(l)\,dl- \frac{(l_{j+1}-l_j)^2}{2} \Phi'(l_j)
+ \cdots
= \int_{l_j}^{l_{j+1}} \Phi(l)\,dl + O(\mu).
\label{eqn-4-6}
\end{equation}
Using \eqref{eqn-4-6} and noting that $l_{j+1}-l_j$ is 
$O(\mu^{1/2})$, we have 
\begin{equation*}
\lambda_n^2 = \lambda_0^2 - \frac{2}{c_1} \int_{l_0}^{l_n}\Phi(l) dl
+ O(\mu^{1/2}).
\end{equation*}
Taking the limit $n\rightarrow\infty$, we find that formally the
stable manifold is given by $\lambda = u(l) + O(\mu^{1/2})$, where
\begin{equation}
u(l)^2 = -\frac{2}{c_1}\int_0^l \Phi(l)\, dl.
\label{eqn-4-7}
\end{equation}
The positive root must be used if $l>0$ and the negative root
if $l<0$. 
It can be verified that this expression for $u(l)$ agrees with
\eqref{eqn-4-4} for the slope at the origin.
Figure \ref{fig-2} plots $u(l)$ and $U(l) = u(l)^2$.

To find the next term in the expansion of the stable manifold,
we use \eqref{eqn-4-5} to get
\begin{equation}
\lambda_n^2 = \lambda_0^2
+ 2\mu^{1/2}\sum_{j=0}^{n-1} \lambda_j \Phi(l_j)
+ \mu \sum_{j=0}^{n-1}\Phi(l_j)^2 + 2\lambda_j^2\Psi(l_j)
+\cdots
\label{eqn-4-8}
\end{equation}
and use \eqref{eqn-4-5} and \eqref{eqn-4-6} to get
\begin{align}
\lambda_j \Phi(l_j) \mu^{1/2} &=
-\frac{1}{c_1} (l_{j+1}-l_j)\Phi(l_j)
+\frac{c_2\lambda_j^2+\Chi(l_j)}{c_1} \Phi(l_j) \mu + \cdots
\nonumber\\
&= -\frac{1}{c_1} \int_{l_j}^{l_{j+1}}\Phi(l)\, dl
+ \frac{1}{c_1}\bigl(c_1^2\lambda_j^2\Phi'(l_j)/2 + (c_2\lambda_j^2
+\Chi(l_j))\Phi(l_j)\bigr)\mu + \cdots
\label{eqn-4-9}
\end{align}
Using \eqref{eqn-4-8} and \eqref{eqn-4-9} and by turning a sum 
into an integral as before, we get the expansion $\lambda
= u(l) + v(l) \mu^{1/2} + O(\mu)$ for the stable manifold ,
where
\begin{equation}
v(l) = -\frac{1}{c_1 u(l)} \int_0^l 
\frac{c_1 u(l)\Phi'(l)}{2} + \frac{c_2u(l)^2+\Chi(l)}{u(l)}\Phi(l)
+\frac{\Phi(l)^2}{2u(l)} + u(l) \Psi(l)\, dl.
\label{eqn-4-10}
\end{equation}
This formula for $v(l)$ agrees with \eqref{eqn-4-4} with regard
to the slope at the origin. This procedure can be repeated to
calculate more terms in the expansion of the stable manifold.

\section{Homoclinic points near resonances}

Assume that the stable manifold of the origin under the map
\eqref{eqn-4-1} is the graph of the function $\lambda = M_0(l)$.
Then $M_0(l)$ must satisfy a functional equation of the form
\begin{equation*}
\lambda_0 + \Phi(l_0)\mu^{1/2} +\cdots
= M_0(l_0-c_1\lambda_0\mu^{1/2}+\cdots).
\end{equation*}
Assume that $M_0(l) = u(l) + v(l) \mu^{1/2} + w(l) \mu$,
where $u(l)$ and $v(l)$ are given by \eqref{eqn-4-7} and \eqref{eqn-4-10}, respectively. Then $w(l)$ must satisfy
the functional equation
\begin{equation}
(1+c_1u'(l_0)\mu^{1/2})w(l_0) = w(l_1) + s_2(l_0,\sqrt{\mu}w(l_0),
\sqrt{\mu})\mu^{1/2},
\label{eqn-5-1}
\end{equation}
where $s_2$ is analytic in its arguments for $0\leq l_0\leq 3\pi/2p$,
$\sqrt{\mu} w$ real and bounded by a large constant, and
$\abs{\sqrt{\mu}} \leq \sqrt{\mu_0}$ for some $\mu_0>0$. In
addition, $s_2(0,0,\sqrt{\mu})=0$. In Lemma \ref{lem-5-1} below , we prove that 
\eqref{eqn-5-1} has a unique $C^1$ solution $w(l)$, with $0\leq l \leq 3\pi/2p$,
for $0<\mu \leq \mu_0$ and some $\mu_0 > 0$. We also prove that $\abs{w(l)}$
and $\abs{w'(l)}$ are bounded by constants which are independent of $\mu$ but which
may depend upon $\mu_0$. 
Therefore, the stable manifold of \eqref{eqn-4-1} is the graph of
$\lambda = u(l) + v(l)\mu^{1/2} + w(l) \mu$ for $0\leq l \leq 3\pi/2p$ and
$0<\mu\leq \mu_0$, where $u(l)$ and $v(l)$ are defined by \eqref{eqn-4-7} and
\eqref{eqn-4-10}.

As discussed in \cite{Zehnder1}, the functional equation \eqref{eqn-5-1} for $w$ 
is obviously a contraction for $0\leq l \leq \pi/p-\epsilon$, $\epsilon > 0$. But the construction of homoclinic points requires the existence of $w$
to be proved over a larger interval as in the lemma below. The proof of the lemma uses
a technique found in \cite{Zehnder1}. Another approach can be found
in \cite{GL1}.

\begin{lem}
For some $\mu_0 > 0$ and any $\mu \in (0,\mu_0]$, there exists a unique $C^1$ 
function $w(l)$ such that \eqref{eqn-5-1} is satisfied for $0\leq l_0\leq 3\pi/2p$,
with $\abs{w(l)}$ and $\abs{w'(l)}$ bounded by constants for $0\leq l \leq 3\pi/2p$.
The constants are independent of $\mu$ but may depend upon $\mu_0$. 
\label{lem-5-1}
\end{lem}
\begin{proof}
We will look for a continuous solution of \eqref{eqn-5-1} that satisfies $w(0) = 0$, $\abs{\exp(-Kl)
w(l)} \leq C^\ast$, and $\lip(\exp(-Kl)w(l))\leq L^\ast$, where $\lip(\cdot)$ is the 
Lipschitz constant. The choice of the positive constants $K$, $C^\ast$, and $L^\ast$ will
be made later in the proof. 

The function $r(l,\lambda, \sqrt{\mu})$ from the first line of \eqref{eqn-4-1}
and the function $s_2(l, \sqrt{\mu}\, w, \sqrt{\mu})$ from \eqref{eqn-5-1} determine
the initial choice of $\mu_0$. Let $w^\ast = \mu \, w$ and consider
$r^\ast(l, w^\ast, \sqrt{\mu}) = r(l, u(l) + \sqrt{\mu}v(l) + w^\ast, \sqrt{\mu})$.
The constant $\mu_0>0$ is chosen so that $r^\ast(l, w^\ast, \sqrt{\mu})$ is analytic
in its arguments 
over the compact domain $D$ where $0\leq l \leq 3\pi/2p$, $w^\ast$ is real and
$\abs{w^\ast} \leq C$, and $\abs{\sqrt{\mu}} \leq \abs{\sqrt{\mu_0}}$. By taking
$\mu_0$ small enough, we can assume $C$ to be as large as we please. It is enough to assume
$C$ to be twice the height of $u(l)$ depicted in Figure \ref{fig-2}, for example.
Now let $w^\ast = \sqrt{\mu}\, w$ and assume the choice of $\mu_0$ to be such that
$s_2(l, w^\ast, \sqrt{\mu})$ is also analytic in its arguments in the compact domain
$D$. We note $r^\ast(0,0,\sqrt{\mu}) = s_2(0,0,\sqrt{\mu}) = 0$.

The proof, which is organized into a number of steps, introduces many constants. The
constants that depend on the domain $D$ will be denoted by subscripting $D$. The constants
that do not depend upon the domain $D$ will be denoted by subscripting $C$. The constants 
that depend upon $D$ are typically upper bounds for the magnitudes of 
derivatives of $r^\ast$ and 
$s_2$ over the domain $D$. The constant $\mu_0$ may be made smaller by some of the steps 
in the proof. But the bounds obtained using the domain $D$ as specified above will of
course apply even if $\mu_0$ is made smaller. All constants introduced in the proof are 
strictly positive.

\begin{enumerate}
\item For $0\leq l_0 \leq 3\pi/2p$, by \eqref{eqn-4-1}
$l_1$ as a function of $l_0$ is given by
\begin{equation}
l_1 = l_0 - c_1 \lambda_0 \mu^{1/2} + (c_2 \lambda_0^2 + \Chi(l_0))\mu 
+ r(l_0, \lambda_0, \sqrt{\mu})\mu^{3/2},
\label{eqn-5-2}
\end{equation}
where $\lambda_0 = u(l_0) + v(l_0)\mu^{1/2} + w(l_0)\mu$. 
It is possible to think of \eqref{eqn-5-2} as defining $l_1$ in terms
of $l_0$ and $w$.
By the assumption
about $\phi(l)$ in Section 3 and \eqref{eqn-4-7},
it follows that $u'(0) > 0$ and
that $u(l) > 0$ for $0<l\leq 3\pi/2p$.
As $u(0) = 0$, there must be a constant $C_1$ such that $c_1u(l) \geq C_1 l$
for $0\leq l \leq 3\pi/2p$. Both the $\mu$ and $\mu^{3/2}$ terms in \eqref{eqn-5-2}
vanish when $l_0=0$ and $w=0$. 
Further, $\abs{w(l_0)} \leq \lip(w)l_0$. 
Thus the  magnitudes of the two terms 
can be upper bounded by $\bigl(C_2 + C_3 \lip(w)\mu\bigr)l_0\mu$ and 
$\bigl(D_1 + D_2 \lip(w)\mu)l_0 \mu^{3/2}$, respectively. Therefore, we may assert
$l_1 \leq \bigl(1-(C_1/2)\mu^{1/2}\bigr)l_0$
for $0\leq l \leq 3\pi/2p$ and $0\leq \mu \leq \mu_0$,
with $\mu_0$ made smaller if necessary. 

\item Let $l^\ast_1$ be obtained using \eqref{eqn-5-2} with 
$l_0$ replaced by $l^\ast_0$ but with the same $w$. Assume $0\leq l_0, l_0^\ast
\leq 3\pi/2p$ and let $C_4 = \lip(c_1u(l))$. Then, as in the previous step, it follows
that $\abs{l_1-l_1^\ast}\leq \bigl(1+2C_4 \mu^{1/2}\bigr)\abs{l_0-l_0^\ast}$
for $0\leq \mu \leq \mu_0$ and $\mu_0$ sufficiently small.

An additional fact about $\abs{l_1-l_1^\ast}$ will be needed.
Assume $0\leq l_0, l_0^\ast\leq \pi/2p$.
By the mean value theorem, $c_1(u(l_0) - u(l_0^\ast)) = c_1u'(\bar{l}) (l_0-l_0^\ast)$,
and by the assumption about $\phi(l)$ in Section 3 and 
\eqref{eqn-4-6}, $u'(l) \geq C_5$ for $0\leq l \leq \pi/2p$. Therefore,
$\abs{l_1 - l_1^\ast} \leq \bigl(1-(C_5/2)\mu^{1/2}\bigr)\abs{l_0-l_0^\ast}$
or simply $\abs{l_1-l_1^\ast} \leq \abs{l_0-l_0^\ast}$ in this situation, for
$0\leq \mu \leq \mu_0$ and $\mu_0$ sufficiently small. 

\item From \eqref{eqn-5-1}, we may obtain the following iteration:
\begin{equation}
\exp(-Kl_0)w_{n+1}(l_0) = F_{w_n}(l_0)\exp(-Kl_1)w_n(l_1)  + 
\frac{\exp(-Kl_0)s_2(l_0, \sqrt{\mu}\,w_n(l_0), \sqrt{\mu})}
{1+c_1u'(l_0)\mu^{1/2}}\, \mu^{1/2},
\label{eqn-5-3}
\end{equation}
where $l_1$ is obtained from $l_0$ using \eqref{eqn-5-2} but with $w$ replaced
by $w_n$, where $w_n$ belongs to the class of functions for $w$ specified at the
beginning of this proof, and where the contraction factor $F_{w_n}$ is given by
\[F_{w_n}(l_0) = \frac{\exp(K(l_1-l_0))}{1+c_1u'(l_0)\mu^{1/2}}.\]
The constant $K$ will be chosen so as to make $F_{w_n}$ a sufficiently strong 
contraction.

First consider $0\leq l_0\leq \pi/2p$. Let the minimum value of $c_1u'(l_0)$ for $l_0$ 
in this range be $C_6$. Since $l_1-l_0\leq - C_1 l_0 \mu^{1/2}/2$ by the first step,
it follows that $F_{w_n}(l_0) \leq \bigl(1-(C_6/2)\mu^{1/2}\bigr)$ for
$0\leq l_0 \leq \pi/2p$, $0\leq \mu \leq \mu_0$, and $\mu_0$ sufficiently small.

Next consider $\pi/2p\leq l_0 \leq 3\pi/2p$. Let the minimum value of $c_1u'(l)$ for
$l_0$ in this range be $-C_7$. Then
\[F_{w_n}(l_0) \leq \frac{\exp(-KC_1\pi\mu^{1/2}/4p)}{1-C_7\mu^{1/2}}.\]
Choose $K$ so that $KC_1\pi/4p \geq 5C_4 + 2 C_7$ and conclude that
$F_{w_n}(l_0) \leq (1-4 C_4 \mu^{1/2})$ for
$\pi/2p\leq l_0 \leq 3\pi/2p$, $0\leq \mu \leq \mu_0$, and $\mu_0$ sufficiently small.

\item By assumption, $\abs{\exp(-Kl)w_n(l)} \leq C^\ast$ for $0\leq l \leq 3\pi/2p$.
Using \eqref{eqn-5-3}, we may upper bound $\abs{\exp(-Kl)w_{n+1}(l)}$ by
$(1-C_8\mu^{1/2})C^\ast + D_3\mu^{1/2}$, where $C_8=\min(C_6/2,4C_4)$ and $D_3$ is an
upper bound of the coefficient of $\mu^{1/2}$ in the last term of \eqref{eqn-5-3}.
By choosing $C^\ast \geq D_3/C_8$, we assert that  $\abs{\exp(-Kl)w_{n+1}(l)}$ is
also upper bounded by $C^\ast$.

\item By assumption,  $\lip\bigl(\exp(-Kl) w_n(l)\bigr) \leq L^\ast$. Let 
\[Q = \exp(-Kl_0)w_{n+1}(l_0) - \exp(-Kl_0^\ast)w_{n+1}(l_0^\ast),\] where
$0\leq l_0^\ast\leq l_0 \leq 3\pi/2p$. We will upper bound $\abs{Q}$.

Using \eqref{eqn-5-3}, both the terms of $Q$ can be replaced by expressions in terms
of $w_n$. The resulting expression for $Q$ equals $A_1B_1 - A_2B_2$ --- where
$A_1 = F_{w_n}(l_0)$, $B_1 = \exp(-k l_1) w_n(l_1)$, $A_2 = F_{w_n}(l_0^\ast)$, 
and  $B_2 = \exp(-k l_1^\ast) w_n(l_1^\ast)$ --- plus another term which equals
the difference of two quantities times $\mu^{1/2}$. This other term will be denoted
by $Q_r \mu^{1/2}$.

To bound $\abs{Q}$, first consider the case $l_0 \geq \pi/2p$. We write
$\abs{A_1B_1 - A_2 B_2} \leq \abs{A_1}\abs{B_1-B_2}+\abs{B_2}\abs{A_1-A_2}$.
By the third step and the assumption about $l_0$, $\abs{A_1} \leq 
\bigl(1-4C_4\mu^{1/2}\bigr)$, and $\abs{B_1-B_2} \leq L^*\abs{l_1-l_1^\ast}
\leq L^\ast \bigl(1+2C_4\mu^{1/2}\bigr)\abs{l_0 - l_0^\ast}$, where the last
inequality follows from the  second step above. 
Therefore, $\abs{A_1(B_1-B_2)}
\leq \bigl(1-C_4\mu^{1/2}\bigr)L^\ast\abs{l_0-l_0^\ast}$. A simple estimate shows
that $\lip(F_{w_n}(l)) \leq C_9 \mu^{1/2}$ for sufficiently small $\mu_0$. 
Therefore, $\abs{B_2(A_1-A_2)} \leq C^\ast C_9 \mu^{1/2} \abs{l_0-l_0^\ast}$.
To upper bound $Q_r$, note that the coefficient of $\mu^{1/2}$ in \eqref{eqn-5-3}
has a Lipschitz constant with respect to $l_0$ that can be bounded as 
$D_4 + D_5 \lip(w_n)\mu^{1/2}$ or $D_4 + D_5 L^\ast \mu^{1/2}$.  Therefore
\[ \abs{Q} \leq \Bigl(\bigl(1-C_4\mu^{1/2}\bigr)L^\ast + C^\ast C_9 \mu^{1/2} 
+ D_4\mu^{1/2} + D_5 L^\ast\mu\Bigr)\abs{l_0-l_0^\ast}.\] 
If $L^\ast\geq 2(C^\ast C_9+D_4)/C_4$, then $\abs{Q} \leq L^\ast \abs{l_0-l_0^\ast}$
for $0\leq \mu \leq \mu_0$ and $\mu_0$ sufficiently small.

Consider the case $0\leq l_0^\ast\leq l_0\leq \pi/2p$. In this case,
the argument is identical to that given in the previous paragraph, except that
the bound on $\abs{A_1}$ must be replaced by $\bigl(1-(C_6/2)\mu^{1/2}\bigr)$
form the third step, and $\abs{B_1-B_2}\leq L^\ast \abs{l_1-l_1^\ast}
\leq L^\ast\abs{l_0-l_0^\ast}$ from the additional fact in the second step. In this
case, if $L^\ast\geq 4(C^\ast C_9+D_4)/C_6$, then $\abs{Q} \leq L^\ast \abs{l_0-l_0^\ast}$
for $0\leq \mu \leq \mu_0$ and $\mu_0$ sufficiently small.

The choice $L^\ast = 2(C^\ast C_9+D_4)\max(1/C_4, 2/C_6)$ implies 
$\lip\bigl(\exp(-Kl)w_{n+1}(l)\bigr)\leq L^\ast$.

\item If $w_n(0) = 0$, then $w_{n+1}(0) = 0$ since $s_2(0,0,\sqrt{\mu})=0$.
This observation together with the choice of $C^\ast$ and $L^\ast$ in the fourth
and fifth steps implies that $w_{n+1}$ belongs to the same class of functions 
as $w_n$. The third step with some other estimates given above implies that
the map $w_n\rightarrow w_{n+1}$ given by \eqref{eqn-5-3} is a contraction
for $0<\mu\leq \mu_0$.
We conclude that there is a unique continuous solution $w$ of \eqref{eqn-5-1} such
that $\abs{\exp(-Kl)w(l)} \leq C^\ast$ and $\lip\bigl(\exp(-Kl)w(l)\bigr)\leq L^\ast$
for $0\leq l \leq 3\pi/2p$. 

\item If $w(l)$ is continuously differentiable its derivative can be easily bounded
in terms of $K$, $C^\ast$, and $L^\ast$. To complete the proof, it suffices to show that
$w(l)$ is continuously differentiable. The standard stable manifold theorem states that
$w(l)$ will be analytic in $l$ in a neighborhood of $l=0$. The stable manifold over
the interval $0\leq l \leq 3\pi/2p$ can be obtained by repeated applying $T_p^{-1}$ to a 
local segment. Therefore $w(l)$ must be continuously differentiable. 
\end{enumerate}

\end{proof}

If $\phi(l_0)$ defined by \eqref{eqn-3-2} satisfies the assumption 
in Section 3, the fixed points of the map $T_p$ given by \eqref{eqn-3-4}
are hyperbolic for $j$ even or for $j$ odd. One of these hyperbolic points was
shifted to the origin in \eqref{eqn-4-1}, and we proved that the stable manifold of
the origin is given by the graph of $\lambda = u(l) + v(l) \mu^{1/2} + O(\mu)$
for $0\leq l\leq 3\pi/2p$. In the $L$-$l$ plane, the stable manifold is the graph
of 
\begin{equation}
L = \Bigl(\frac{p}{q}\Bigr)^{1/3}  + \frac{\chi(j\pi/p)}{c_1}\, \mu
+ u(l-j\pi/p)\mu^{1/2} + v(l-j\pi/p)\mu + O(\mu^{3/2})
\label{eqn-5-4}
\end{equation}
for $0\leq l-j\pi/p+O(\mu)\leq 3\pi/2p$,
where $\chi$ is given by \eqref{eqn-3-2},  $u$ is given by \eqref{eqn-4-7},
and $v$ is given by \eqref{eqn-4-10}.

The choice of the Poincar\'e section as either $g=0$ or $g=\pi$ is yet to be made.
To facilitate the construction of homoclinic points, it is also useful to pick $j$ in \eqref{eqn-3-4}
carefully. There are four cases.
\begin{itemize}
\item If $p$ is odd, then $g=0$ is chosen as the Poincar\'e section. If 
$\phi'(0) > 0$, then $j=-1$.
\item If $p$ is odd and $\phi'(0) < 0$, then $j=(p-1)$.
\item If $p$ is even, first try $g=0$ as the Poincar\'e section. If 
$\phi'(0) > 0$, then $j=-1$.
\item If $p$ is even and $\phi'(0) < 0$ with $g=0$ as the Poincar\'e section,
choose the Poincar\'e section $g=\pi$ and $j=-1$.
\end{itemize}

\begin{thm}
Let $p$ and $q$ be  relatively prime positive integers and let $p/q\neq 1/1$. 
Assume that $e$ lies in the interval $[e_{min}, e_{max}]$ defined at the beginning of Section 2. Let $\phi(l_0)$ be defined by \eqref{eqn-3-2} with $L=(p/q)^{1/3}$ and $G = (p/q)^{1/3}(1-e^2)$. Assume 
that $\phi(l_0)$ satisfies Assumption A of 
Section 3,
namely, for $l_0\in [0,2\pi/p)$, $\phi(l_0) = 0$ only if $l_0=0$ or $l_0=\pi/p$
and $\phi'(l_0)\neq 0$ at those two points. Identify the Poincar\'e section
for the flow of the Hamiltonian \eqref{eqn-1-2} of the restricted three-body problem
with a region of the 
$L$-$l$ plane by using $H = -(p/q)^{-2/3}/2 - (p/q)^{1/3}(1-e^2)$ and by choosing
$g=0$ or $g=\pi$ as indicated above. Then the $p$th return map $T_p$ given by
\eqref{eqn-2-2} has a homoclinic point on this Poincar\'e section at $(l_h,L_h)$,
where $L_h = (p/q)^{1/3} + u(\pi/p)\mu^{1/2}+O(\mu)$ and  $l_h=0$ or $l_h=\pi$, for $0<\mu\leq \mu_0$
and $\mu_0$ sufficiently small.
\label{thm-5-1}
\end{thm}
\begin{proof}
We give a proof for the second case listed above. The other cases are treated similarly.
In this case, $p$ is odd, the Poincar\'e section is $g=0$, and $j=(p-1)$.
By Lemma \ref{lem-5-1} the representation of the stable manifold given by \eqref{eqn-5-4} is valid
for $0\leq l-(p-1)\pi/p + O(\mu) \leq 3\pi/2p$, and therefore the stable manifold
crosses the line $l=\pi$. By Lemma \ref{lem-2-1}, this stable manifold can be reflected
about the line $l=\pi$ to obtain an unstable manifold.  Thus we find a homoclinic 
point with $l=\pi$. 
\end{proof}

The homoclinic point constructed in Theorem \ref{thm-5-1} can be mapped using
the first return map $T_1$ to obtain a ring of $p$ homoclinic points with 
$L_h > (p/q)^{1/3}$. In Section 4, we constructed the stable manifold of the fixed
point of \eqref{eqn-4-1} at the origin over $0\leq l_0 \leq 3\pi/2p$. A similar
construction applies over the interval $-3\pi/2p\leq l_0\leq 0$. That construction
can be used to find a ring of $p$ homoclinic points with $L_h < (p/q)^{1/3}$.

\section{Verification of the condition on $\phi(l_0)$}
Equations \eqref{eqn-1-3} and \eqref{eqn-1-4} define $\Omega$ as a function of
$L,l,G,g$. If the angles $l$ and $g$ are replaced by $-l$ and $-g$, $\Omega$ is unchanged.
Therefore, $\Omega$ can be Fourier expanded as
\begin{equation}
\Omega = \sum_{m,n} c_{mn} \cos(ml + ng), 
\label{eqn-6-1}
\end{equation}
where $m$ can be any non-negative integer and $n$ can be any integer.
The coefficients $c_{mn}$ are functions of $L$ and $G$. We use $L=(p/q)^{1/3}$
and $G=(p/q)^{1/3}(1-e^2)^{1/2}$. The Fourier expansion is valid if
$1+r^2-2r\cos(\theta) > 0$ and is therefore valid if $e$ is sufficiently small
for $p/q\neq 1/1$. The coefficients $c_{mn}$ can be expanded as power series in $e$,
and it is possible to determine the precise radius of convergence of these series. For our purposes, it suffices to note that all these series converge in some neighborhood
of $e=0$. If $\Omega$ is differentiated with respect to $L$ or $G$ the condition
$e>0$ has to be imposed. However, the partial derivatives of $\Omega$ with respect to $l$
and $g$ and $\Omega$ itself are analytic in a neighborhood of $e=0$.

\begin{lem}
If the Fourier coefficients $c_{mn}$ of \eqref{eqn-6-1} are expanded in powers of $e$,
the lowest power of $e$ with a possibly nonzero coefficient is $e^{\abs{m-n}}$.
\label{lem-6-1}
\end{lem}
\begin{proof}
The quantity $\Omega$ can be expanded as $\sum_{m=0}^\infty c_m(r) \cos(m\theta)$.
First consider $e=0$. Then $r = (p/q)^{1/3}$ and $\theta = \nu + g = l+g$.
Therefore the only nonzero terms in the expansion \eqref{eqn-6-1} occur when
$m=n$. If $e\neq 0$, then $r$ depends upon $l$ and $\nu$ is no longer equal to
$l$. This lemma follows when the dependence of $r$ and $\nu$ on $l$ is taken into
account. The way to do this can be found on pages 44 and 170 of \cite{Plummer1}
or on page 35 of \cite{Morbidelli1}.
\end{proof}

By \eqref{eqn-3-2}, $\phi(l_0) = \int_0^{2\pi p} \Omega_l\, dt$, where
$\Omega_l$ must be evaluated at $l=l_0+qt/p$, $g=g_0-t$, $L=(p/q)^{1/3}$,
and $G=(p/q)^{1/3}(1-e^2)^{1/2}$; $g_0=0$ or $g_0=\pi$ depending upon the choice
of the Poincar\'e section. Using \eqref{eqn-6-1}, we get
\[ \phi(l_0) = -2\pi p^2\sum_{m=1}^\infty m\, c_{mp,mq} \sin(kpl_0 + kq g_0).\]
Let $c_{p,q} = c^\ast(p,q) e^{\abs{p-q}} + O(e^{\abs{p-q}+1})$. Then by Lemma 
\ref{lem-6-1},
$\phi(l_0) = \pm 2\pi p^2 c^\ast(p,q) \sin(pl_0) e^{\abs{p-q}} + \cdots$.
Thus the assumption about $\phi(l_0)$ in Section 3 and
in Theorem \ref{thm-5-1} will be verified for $e>0$ and $e$ small if we can show
that $c^\ast(p,q) \neq 0$. An expression for $c^\ast(p,q)$ can be obtained from
the discussion of the quantity
\[ C(e,p,q) = - \frac{6\pi q^{4/3}}{p^{1/3}} \int_0^{2\pi p} \Omega_{ll}\, dt\]
given in \cite{Viswanath1}. If $p<q$, for example, 
\[c^\ast(p,q) = - \frac{(-1)^{q-p} q^{2/3}}{6.2^{q-p} \pi p^{8/3}}
\Biggl(\sum_{k=0}^{q-p} \binom{D+q}{k} \frac{p^{q-p-k}}{(q-p-k)!}\Biggr)
(\alpha b_q(\alpha))\]
evaluated at $\alpha = (p/q)^{1/3}$; above $D$ stands for the differential
operator $\alpha \frac{d}{d\alpha}$ and $b_n(\alpha)$ are defined by
the expansion 
\[(1+\alpha^2-2\alpha\cos \theta) = \frac{1}{2}\sum_{n=-\infty}^\infty b_n(\alpha)\exp(in\theta).\]
The $b_n(\alpha)$ are hypergeometric functions whose series converge for
$\abs{\alpha} < 1$.  The value of $c^\ast(p,q)$ can be obtained using the expression
given above or by other means.

\begin{figure}
\begin{center}
\includegraphics[scale=0.5]{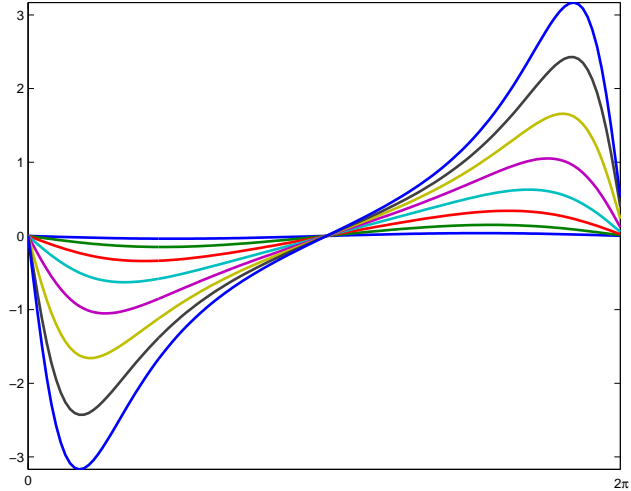}
\end{center}
\caption[xyz]{Plots of $\phi(l)$ against $l$ with $q/p=3/1$ and with
$e$ ranging from $0.1$ to $0.8$.} 
\label{fig-3}
\end{figure}

In Figure \ref{fig-3}, we have plotted $\phi(l)$ with $q/p=3/1$. From that figure,
it is clear that the assumption about $\phi(l_0)$ is valid for even large
values of $e$.

\section{Resonance boundaries}
Discussion of the averaged circular restricted problem and its use
in sketching the boundaries of $q/p$ Kirkwood gaps in the $a$-$e$ plane
can be found in \cite{DM1}, \cite{HL1}, \cite{Lemaitre1}, 
\cite{Morbidelli1}, and \cite{Yoshikawa1}. The basic procedure
is to average the Hamiltonian \eqref{eqn-1-2} by retaining only the
terms in the Fourier expansion \eqref{eqn-6-1} of $\Omega$ with
$ml+ng = k(pl-qg)$ for some integer $k$.
These are the resonant terms.
In some instances such
as \cite{Yoshikawa1}, certain additional terms are added to model the
effect of the eccentricity of Jupiter and the secular variation of
its elements.

This averaged Hamiltonian has 1 degree of freedom. Its fixed points
and separatrices are used to approximate the boundaries of resonance. 
For small values of the asteroid eccentricity $e$,
some of the nonresonant 
terms dropped during averaging have larger coefficients
than any of the terms retained during averaging. To some extent 
the influence of the nonresonant terms is captured by the formal
change of variables used to average the Hamiltonian, but this
change of variables is often not taken into account. Even if it
is, the averaging will not be valid at small values of $e$.

The return map \eqref{eqn-3-1} corresponds to the unaveraged
circular restricted problem. Only the resonant terms of $\Omega$
contribute to $\phi(l)$ defined by \eqref{eqn-3-2}, but the
nonresonant terms contribute to both $\psi(l)$ and $\chi(l)$.
As $e$ approaches $0$, the magnitude of $\phi(l)$ becomes much
smaller than that of the other two functions in \eqref{eqn-3-2}.
Thus for fixed $\mu$ and small $e$ the periodic points 
\eqref{eqn-3-4} will not exist.

These periodic points, when they exist, correspond to the fixed
points of the averaged Hamiltonian. The separatrices of the
averaged Hamiltonian corresponds to the stable manifolds discussed
in Sections 4 and 5. If for certain values of $\mu,e,p,q$  the
unaveraged circular restricted problem does not have the periodic
points given by \eqref{eqn-3-4}, the fixed points and separatrices
of the averaged Hamiltonian must be treated as artifacts of 
the averaging procedure.

\begin{table}
\begin{center}
\begin{tabular}{c|c|c|c|c|c|c}
$q/p$ & $3/2$ & $2/1$ & $7/3$ & $5/2$ & $3/1$ & $4/1$ \\\hline
$e$ & $.10$ & $.13$ & $.08$ & $.08$ &  $.07$ & $.09$ 
\end{tabular}
\end{center}
\caption[xyz]{ This table gives the values of $e$ for certain
$q/p$ below which the periodic points of Section 3 fail to exist
for $\mu = 10^{-3}$. Boundaries of resonance obtained using the
averaged circular restricted problem will not be valid below these
values.}
\label{table-1}
\end{table}

In Table \ref{table-1}, we have given the minimum values of $e$
required for the periodic points \eqref{eqn-3-4} to exist
for some of the commonly studied resonances in the asteroid belt.
The value of $\mu$ used is close to that of Jupiter. The boundaries
of resonance obtained by averaging can be valid only above these
values of the asteroid eccentricity.

\section{Asymmetric librations}

If the assumption about $\phi(l)$ in Section 3 holds, the circular
restricted problem has two resonant periodic solutions for
$\mu$ small. One of these is of elliptic type and therefore there will
be solutions that librate around the periodic points that correspond
to it in the Poincar\'e section given by the $l$-$L$ plane.  These
are symmetric librations.

\begin{figure}
\begin{center}
\includegraphics[scale=0.4]{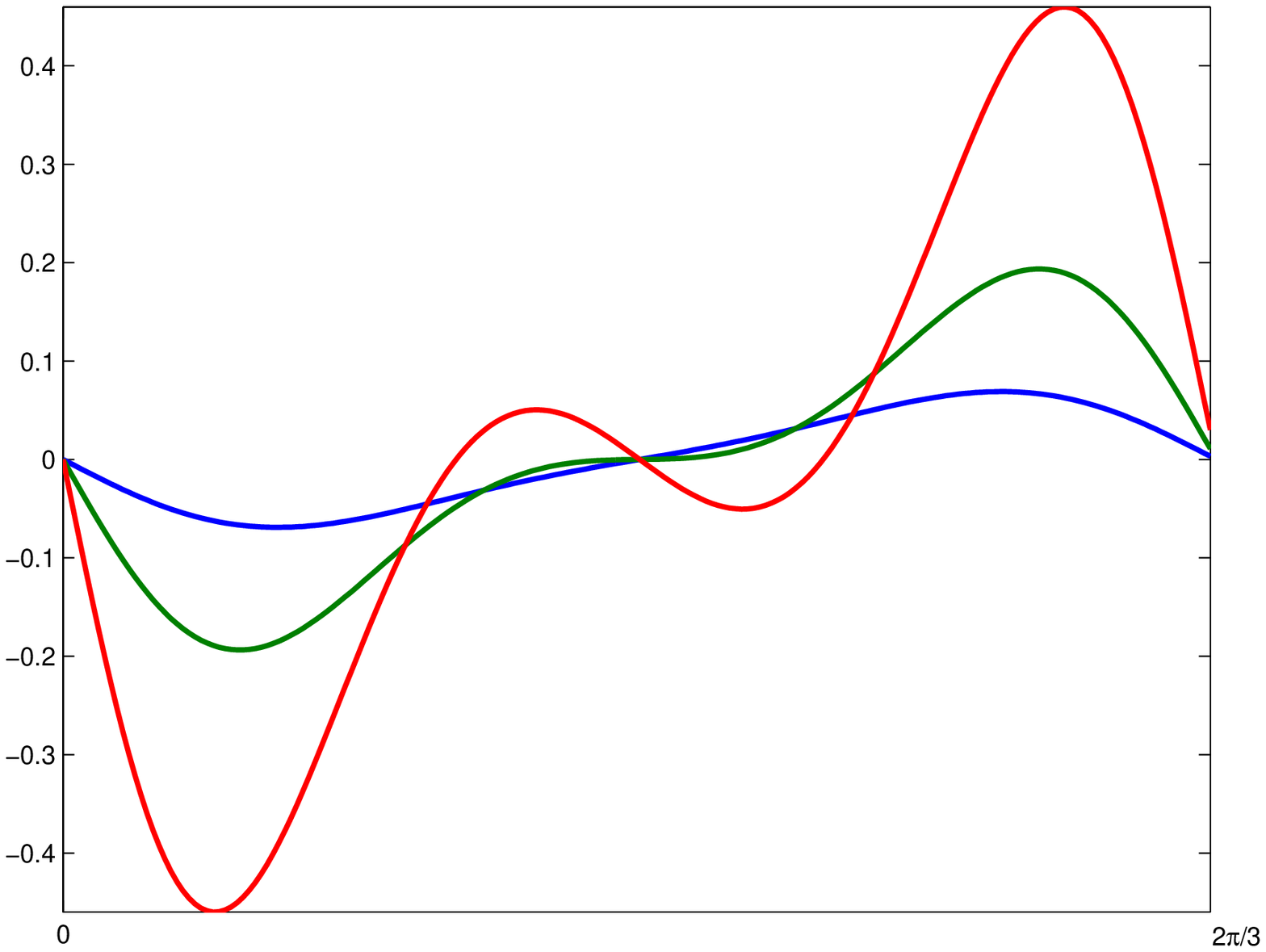}\hspace{.25cm}
\includegraphics[scale=0.4]{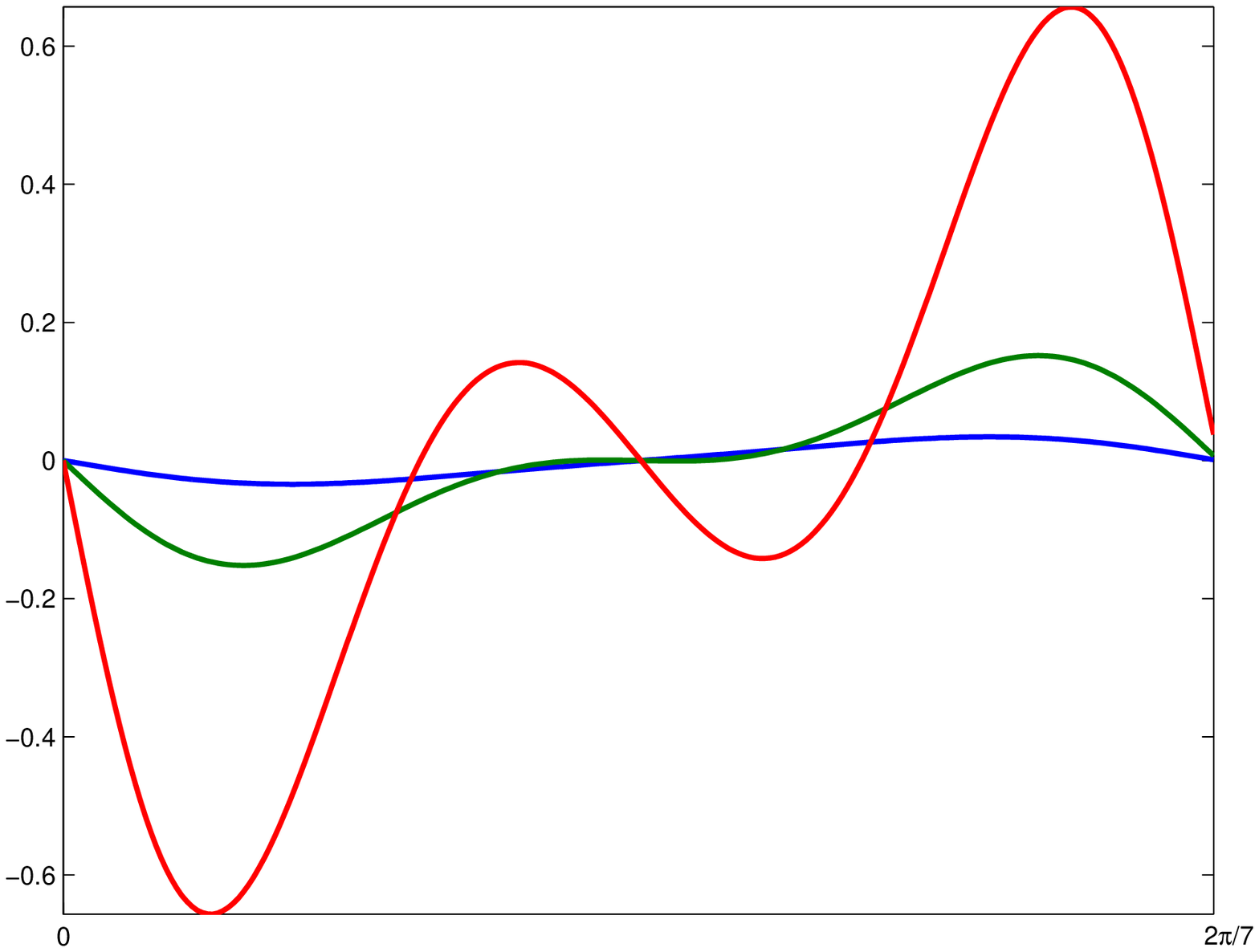}
\end{center}
\caption[xyz]{The two figures show plots of $\phi(l)$ vs. $l$
for $l\in[0,2\pi/p)$ for $q/p=1/3$ and $q/p=1/7$, respectively.
The values of $e$ are $0.08$, $0.12$ and $0.16$ in the left plot,
and $.30$, $.36$ and $.40$ in the right plot.}
\label{fig-4}
\end{figure}

As shown in Figure \ref{fig-4}, the assumption about $\phi(l)$
in Section 3 can fail for some exterior resonances as $e$
increases. In both the plots shown in that figure, the periodic
point with $l\approx \pi/p$ is initially of elliptic type as
the slope of $\phi(l)$ is positive, and when $e$ increases, it 
undergoes a pitchfork bifurcation and becomes a point of hyperbolic type. Elliptic points appear on the
Poincar\'e section at values of $l$ that are not $O(\mu)$ close
to any integral multiple of $\pi/p$. The librations around these
points are termed asymmetric.

\begin{table}
\begin{center}
\begin{tabular}{c|c|c|c|c|c|c}
$q/p$ & $1/7$ & $1/6$ & $1/5$ & $1/4$ & $1/3$ & $1/2$ \\\hline
$e$ & $.365900$ & $.320133$ & $.265532$ & $.199749$ &  $.121094$ & 
$.036083$ 
\end{tabular}
\end{center}
\label{table-2}
\caption[xyz]{This table gives the values of $e$ above which 
asymmetric librations exist in the unaveraged circular restricted
problem for certain exterior resonances.}
\end{table}

A study of asymmetric librations in the exterior resonances using
averaged equations can be found in \cite{Beauge1}. In Table
\ref{table-2}, we have listed values of $e$ above which asymmetric
librations occur in the unaveraged equations for some exterior resonances. If the averaged equations imply the existence of 
asymmetric librations below these values of $e$, those must
be considered artifacts of averaging. However, if $e$ exceeds these
values, asymmetric librations will be found for $\mu$ small enough.

\section{Acknowledgments}
The author thanks Profs. J.E. Marsden and D. Scheeres for helpful
discussions.

%\bibliography{references}
%\bibliographystyle{plain}

\end{document}